\documentclass[11pt]{amsart}

 \newtheorem{thm}{Theorem}[section]
 \newtheorem{con}[thm]{Conjecture}

 \theoremstyle{definition}
 
 \theoremstyle{remark}
 
 \numberwithin{equation}{section}
 \newcommand{\eps}{\varepsilon}
 \newcommand{\To}{\longrightarrow}

\begin{document}
 
 \title{Sharp bounds for eigenvalues of  triangles}
 
 \author{Bart{\l}omiej Siudeja}

\address{Department of Mathematics, Purdue University, West Lafayette, Indiana
47906}

\email{siudeja@math.purdue.edu}


\thanks{Supported in part by NSF Grant \# 9700585-DMS}


\begin{abstract}
We prove that the first eigenvalue of the Dirichlet Laplacian for a triangle in the plane is bounded above 
by $\pi^2 L^2\over 9A^2$, where $L$ is the perimeter and $A$ is the area of this triangle. We show
that the \mbox{constant 9} is optimal and that
the optimal constant for the lower bound of the same form is $16$.
This gives a positive answer to a conjecture made in  \cite{F}.
\end{abstract}

\maketitle
 
 \section{Introduction}

The purpose of this paper is to prove the following theorem. 
 \begin{thm}\label{main}
  Let $T$ be a triangle in a plane of area $A$ and perimeter $L$. Then the first eigenvalue $\lambda_T$ of
  the Dirichlet Laplacian on $T$ satisfies
  \begin{gather}
     {\pi^2 L^2\over 16A^2}\leq\lambda_T\leq {\pi^2 L^2\over 9A^2}. 
  \end{gather}
 The constants $9$ and $16$ are optimal.
\end{thm}

The lower bound was proved in a more general context in \cite{M}. In Section 6 we show that for ``tall"  isosceles 
triangles there is an asymptotic equality in the lower bound. Hence it is impossible to decrease the constant $16$.

The upper bound was recently stated as a conjecture in \cite{F} and numerical evidence for 
its validity are given in \cite{AF}. 
Bounds of this form but with different constants have been the subject of many papers
in the literature.
The eigenvalue of any doubly-connected domain
is bounded above by the same fraction but with the constant $4$,  see \cite{Po} and remarks in \cite{O}. 
There is also a sharper upper bound due to Freitas (\cite{F}) which is not of this form 
but it seems that in the worst case (``tall'' isosceles triangle) it gives the constant $6$, in the best (equilateral) $9$.
It is worth noting that the constant $9$ may not be improved since equilateral triangles give equality in the upper bound of 
the Theorem \ref{main}.

The spectral properties of a Dirichlet Laplacian on an arbitrary planar domain are important both in physics 
and in mathematics.  Unfortunately,
it is almost impossible to find the exact spectrum even for some simple classes of domains. Except for rectangles,
balls and annuli, not much can be said in general. In the case of triangles the full spectrum is known only for equilateral  and right triangles
with smallest angles $\pi/4$ or $\pi/6$. For more information about these
we refer the reader to \cite{Mc}. For all other triangles, the best we can hope to do is to give bounds for the 
eigenvalues, such as those given above. 

Even though Theorem \ref{main} gives sharp bounds in the sense that the constants are the best possible given the form
of the bound, there is certainly room for improvements. In fact, sharper lower bounds are already known, see
\cite{F}. One of these bounds is good for both equilateral and ``tall`` triangles. It gives the constant $9$ for the
first and $16$ for the second. 
The upper bound good in both cases is still unknown. To the best of our knowledge
it is also not clear what is the correct bound for the isosceles triangle with base almost equal to the half of the
perimeter, but we think it should be $16$.

By comparing our numerical results with the numerical studies contained in \cite{AF}, Section 5.1, we conjecture that
\begin{con}  Let $T$ be a triangle in a plane of area $A$ and perimeter $L$. Then the first eigenvalue $\lambda_T$ of
  the Dirichlet Laplacian on $T$ satisfies
\begin{gather}
  {\pi^2L^2\over16A^2}+{7\sqrt3\pi^2\over12A}\leq\lambda_T\leq {\pi^2L^2\over12A^2}+{\sqrt3\pi^2\over3A}.
\end{gather}
\end{con}

Here both bounds are of the form 
$$E_3(L,A,\theta)={4\pi^2\over \sqrt3A}+\theta
{L^2-12\sqrt3A \over A^2}$$ considered in \cite{AF}.

The lower bound, with $\theta=\pi^2/16$, is the best bound we can expect
given this particular form. Indeed, this is the only bound which is sharper then the lower bound of Theorem \ref{main} and which might
be true for ``tall'' triangles.
The upper bound from our main result is also of this form, but with $\theta=\pi^2/9$. Hence the conjectured upper bound
is sharper ($\theta=\pi^2/12$), and it is the best in the sense that the bound with $\theta=\pi^2/13$ is not valid.
But, since only the constant $16$ can give a good upper bound for ``tall'' triangles, it is not possible to find a bound of the form $E_3$ which
is good for both equilateral and ``tall'' triangles.

Our proof of the upper bound from the Theorem \ref{main} contains two main parts. The first deals with an ``almost equilateral'' triangles. That is,  with the triangles for which
the longest side is comparable to the  shortest side. For these our strategy is to find a suitable test function $\psi$. That is, we try to find a function which 
is $0$ on the boundary of the triangle $T$ and apply the Rayleigh quotient to get the upper bound for 
$\lambda_T$.  We get 
\begin{gather}
  \lambda_T\leq {\int_T |\nabla\psi|^2\over \int_T \psi^2}.
\end{gather}
This part of the proof is contained in Sections 2 to 5. Included in Section 2 are also some 
preliminary results. 

The second part of the proof, contained in Section 6, deals with ``tall'' triangles. These can be  approximated by a circular sections for which the 
eigenvalues can be found explicitly.

\section{Eigenfunctions and notation}

An arbitrary triangle $T'$ can be rotated and rescaled to obtain a triangle $T$ with vertices $(0,0)$, $(1,0)$ 
and $(a,b)$. This, together with the fact that the bound in the main theorem is invariant under translations, rotations and scaling
allow us to restrict our attention to the triangles with such vertices. We can also assume that the side 
contained in the $x$-axis is the shortest. Hence we have that  $$a^2+b^2\geq1\,\,\, \text{and}\,\,\,  a\leq1/2$$
for our triangles.  We will denote the length of the other two sides by 
$M$ and $N$, with $N$ denoting the longest.

We start with the first eigenfunction of an equilateral triangle, and we will proceed as in \cite{F}. Such function is given by 
\begin{gather}
  f(x,y)=\sin\left(4\pi y\over \sqrt3\right)-\sin\left[2\pi\left( x+{y\over\sqrt3} \right)\right]+
          \sin\left[2\pi\left( x-{y\over\sqrt3} \right)  \right].
\end{gather}

We can compose $f$ with a linear transformation to obtain a function $\phi$ which is equal to $0$ on the boundary of $T$. Namely consider 
\begin{gather}
  \begin{split}
  \phi(x,y)&=f\left(x-{a-1/2\over b}y,{\sqrt3\over 2b}y\right)\\&=
  \sin\left(2\pi y\over b\right)-\sin\left[2\pi\left( x+{(1-a)y\over b} \right)\right]+
          \sin\left[2\pi\left( x-{a y\over b} \right)  \right].
  \end{split}
\end{gather}

This function was used in \cite{F} to obtain the upper bound from the Rayleigh quotient. Since the function $f$ is the first eigenfunction 
of the Dirichlet Laplacian on an equilateral triangle and its eigenvalue gives equal sign in the main bound, it is reasonable
to expect that by taking any linear transformation we can only decrease the constant $9$ in the Theorem \ref{main}. 

Hence we want to  find another eigenfunction of some other triangle. We will use the eigenfunctions of the equilateral triangle to
find a test function for the right triangle with angles $\pi/3$ and $\pi/6$. In the recent paper \cite{Mc} the author constructs two families
of eigenfunctions of the equilateral triangle. The antisymmetric mode has the property that it is $0$ on the altitude.
Thus, such a function is also the eigenfunction for the right triangle. We can then take the antisymmetric eigenfunction 
corresponding to the smallest eigenvalue as our test function. A calculation leads to the following function
\begin{gather}
  \begin{split}
  g(x,y)&=\sin\left( \sqrt3\pi y \right)\sin\left( \pi x/3 \right)
  \\&+
  \sin\left( \sqrt3\pi y/3 \right)\sin\left( 5\pi x/3 \right)
  \\&+
  \sin\left( 2\sqrt3\pi y/3 \right)\sin\left( 4\pi x/3 \right).
  \end{split}
\end{gather}

This function, as can be easily checked, is in fact the eigenfunction of the Dirichlet Laplacian on the triangle with the vertices
$(0,0)$, $(1,0)$ and $(0,\sqrt3)$. The corresponding eigenvalue gives the better bound than the one in the Theorem \ref{main}, the
constant is about $9.6$.  Therefore, a linear transformation of this function should give a correct bound at least for the
neighborhood of the point 
$(0,\sqrt3)$. By applying a suitable linear transformation we get the second test function
\begin{gather}
  \begin{split}
    \varphi_1(x,y)&=g\left( x-{ay\over b},{\sqrt3 y\over b} \right)
   \\&=
    \sin\left( 3\pi y\over b \right)\sin\left[\frac\pi3\left( x-{ay\over b} \right) \right]
    \\&+
    \sin\left(\pi y\over b \right)\sin\left[ \frac{5\pi}3\left(x-{ay\over b}  \right) \right]
  \\&+
  \sin\left( 2\pi y\over b \right)\sin\left[ \frac{4\pi}3\left(x-{ay\over b}  \right) \right].
  \end{split}
\end{gather}

Similarly, we can obtain the last two test functions. One will be a linear transformation of the eigenfunction of the triangle 
with vertices $(0,0)$, $(1,0)$ and $(1,\sqrt3)$. The other a linear transformation of the eigenfunction of the 
triangle with the vertices $(0,0)$, $(1,0)$ and $(0,1/\sqrt3)$. We get:
\begin{gather}
  \begin{split}
    \varphi_2(x,y)&=
    \sin\left( 3\pi y\over b \right)\sin\left[\frac\pi3\left( 1-x+{(a-1)y\over b} \right) \right]
    \\&+
    \sin\left(\pi y\over b \right)\sin\left[ \frac{5\pi}3\left(1-x+{(a-1)y\over b}  \right) \right]
  \\&+
  \sin\left( 2\pi y\over b \right)\sin\left[ \frac{4\pi}3\left(1-x+{(a-1)y\over b}  \right) \right].
  \\
    \varphi_3(x,y)&=
    \sin\left( 5\pi y\over 3b \right)\sin\left[\pi\left( x-{ay\over b} \right) \right]
    \\&+
    \sin\left(4\pi y\over 3b \right)\sin\left[2\pi\left(x-{ay\over b}  \right) \right]
  \\&+
  \sin\left( \pi y\over 3b \right)\sin\left[3\pi\left(x-{ay\over b}  \right) \right].
  \end{split}
\end{gather}

Now we can take a linear combination of these test functions.  That is, we consider 
\begin{equation}
  \psi(x,y)=\alpha\varphi_1(x,y)+\beta\varphi_2(x,y)+\gamma\varphi_3(x,y)+\eps\phi(x,y),
\end{equation}
and we can calculate the Rayleigh quotient for this function.  After optimizing over all possible values of 
$\alpha$, $\beta$, $\gamma$ and $\eps$, this will give an appropriate bound for the first eigenvalue.
To prove Theorem \ref{main} we have to check that
\begin{gather}
  \lambda_T\leq {\int_T|\nabla\psi|^2\over\int_T\psi^2}\leq {\pi^2 L^2\over 9A^2},
\end{gather}
for some $\alpha$, $\beta$, $\gamma$ and $\eps$ (possibly depending on $T$).
The last  inequality is equivalent to
\begin{equation}\label{in1}
 9A^2\int_T|\nabla\psi|^2\leq \pi^2 L^2\int_T\psi^2. 
\end{equation}

Since the function $\psi$ is given explicitly and is a trigonometric function, it is possible to find the exact values of
this integrals but calculations are very cumbersome. For this reason we will do the
long calculations  in Mathe\-matica. However, we wish to emphasize the fact that all the calculations are done symbolically.

By our assumptions we have $L=1+\sqrt{a^2+b^2}+\sqrt{(a-1)^2+b^2}$ and $A=b/2$. 
As a result of running Mathematica we get that to prove the inequality (\ref{in1}) we have to find 
$\alpha$, $\beta$, $\gamma$ and $\eps$ such that the inequality  
\begin{small}
\begin{gather}\label{horrible}
  \begin{split}
&0\geq8041366333\times
\\&\;\;\Big\{
\left( -1594323 - 1792090a + 531441(a^2+b^2) + 
     201600\left( 3 + a^2 + b^2 \right) {\pi }^2 \right)\alpha^2 
\\&\;\;\;+ 
  \left( -2854972 + 729208a + 531441(a^2+b^2) + 
     201600\left( 3 + (a-1)^2 + b^2 \right) {\pi }^2
     \right)\beta^2 
\\&\;\;\;+  
   \left( 531441 - 1792090a - 1594323(a^2+b^2) + 
     201600\left( 1 + 3a^2 + 3b^2 \right) {\pi }^2 \right) 
     {\gamma }^2\Big\}
\\&\;\;+ 
  5558192409369600\left( 1 -a + (a^2 + b^2) \right) 
   \pi^2\eps^2
\\&\;\;+67672797192\times
\\&\;\;\;\Big\{
\left( 729{\sqrt{3}}
      \left( 454 -128a + 339(a^2+b^2) \right)  + 
     24640\left( 4 -8a + 9(a^2+b^2) \right) \pi 
     \right) \gamma \epsilon   
\\&\;\;\;\;+
     \left( 729{\sqrt{3}}
         \left( 665 - 780a + 454(a^2+b^2) \right)  + 
        24640\left( 5 + 4(a^2+b^2) \right) \pi  \right)\beta \epsilon 
\\&\;\;\;\;+
    \left( 729{\sqrt{3}}
         \left( 339 -128a + 454(a^2+b^2) \right) 
         + 24640\left( 9 -8a + 4(a^2+b^2) \right) \pi  \right) \alpha\epsilon\Big\}
\\&\;\;+\left( 1990033124626008a + 
        2553294638054160{\sqrt{3}}
         \left( 3 - 2a + 3(a^2+b^2) \right) \pi  \right) \alpha \gamma
\\&\;\;+1151172000\left( 35341051 - 26756686a + 32479596(a^2+b^2) \right) \beta \gamma
\\&\;\;+ 189\Big\{ 819452341268271 + 
        -73323642839420a + 73323642839420(a^2+b^2) 
\\&\;\;\;\;\;\;\;\;\;\;\;\;\;\;- 
        79935610875120{\sqrt{3}}\pi  + 
        24336134222400{\sqrt{3}}
        \left(-a + a^2 + b^2 \right) \pi  \Big\} 
      \alpha \beta 
\\&\;\;
    - 9{\left( 1 + {\sqrt{{\left( -1 + a \right) }^2 + b^2}} + 
       {\sqrt{a^2 + b^2}} \right) }^2
   \Big\{ 444001222376712{\sqrt{3}}\left( \alpha  + \gamma  \right) 
      \epsilon  
\\&\;\;\;\;\;\;\;\;- 1629547920\pi 
      \left( {\sqrt{3}}\alpha 
         \left( 4251\beta  - 99484\gamma  \right)  - 
        113696\left( \alpha  + \beta  + \gamma  \right) \epsilon  \right)
\\&\;\;\;\;\;\;\;\;+ 51464744531200{\pi }^2
      \left( {\alpha }^2 + {\beta }^2 + {\gamma }^2 + 2{\epsilon }^2
        \right)
\\&\;\;\;\;\;\;\;\;+ 3\beta \left( 346474423262177\alpha  + 
        85272\left( 3297684500\gamma  + 
        1735627257{\sqrt{3}}\epsilon  \right)  \right)  \Big\}
  \end{split}
\end{gather}
\end{small}
is valid.

This expression clearly shows that it would be very difficult to do the calculations 
by hand.  Notice that this expression depends only on $b^2$ and $a$. Also note that the ``building" blocks 
for the expression are exactly equal to the length of the sides of the triangle $T$. Hence we make the substitution
$M^2=a^2+b^2$ and $N^2=(a-1)^2+b^2$. As a result we get a polynomial of degree $2$ in $M$  and $N$, where $N\geq M\geq
1$. For further simplification (and to improve our chances of finding the appropriate coefficients) we divide all triangles
into 4 classes. Each class will be handled in a separate section whose number corresponds to the case:
\begin{enumerate}
  \item[3)] The triangles with $N\geq2$ and $M\leq15$,
  \item[4)] $1\leq N\leq2$ and $(N+1)/2\leq M\leq 2$,
  \item[5)] $1\leq N\leq2$ and $1\leq M\leq (N+1)/2$,
  \item[6)] $M\geq 15$.
\end{enumerate}
The method used to handle the last case  will be totally different than the previous ones. 
\section{Case: $N\geq2$ and $M\leq15$}
Let us take $\eps=\beta=0$, $\alpha=1$ and $\gamma=-1/6$. This simplifies
(\ref{horrible}) to
\begin{gather}
  \begin{split}
0\geq &P(M,N)= -90851035780 - 16374894040M^2 + 33929984593N^2 
\\&- 
  272432160{\sqrt{3}}\left( 10 - 8M + 10M^2 - 8N - 8MN + 
     3N^2 \right) \pi  
\\&+ 28828800
   \left( 689 - 148M + 199M^2 - 148N - 148MN - 74N^2 \right) 
   {\pi }^2.
  \end{split}
\end{gather}

To show this inequality we first find all the critical points of the right side and 
later check the values on the boundary. Both $\partial_M P$ and $\partial_N P$ are 
linear with respect to $M$ and $N$, therefore we have exactly $1$ critical point with 
$N\approx-42.2$. Hence it is enough to check this inequality on the boundary.

The boundary conditions are given by $M=N$, $M=15$, $N=2$ and $M=N-1$. For each of these
$P$ is a quadratic equation and we just have to check that the roots are outside of 
the bounds for $M$ or $N$ and that the  inequality is true at the endpoints. We have 
\begin{itemize}
  \item $P(M,M)=0$ for $M\approx1.6$ and $M\approx15.15$; $P(2,2)<0$ and $P(15,15)<0$,
  \item $P(15,N)=0$ for $N\approx14.97$ and $N\approx42.5$; $P(15,15)<0$ and $P(15,16)<0$,
  \item $P(M,2)=0$ for $M\approx0.96$ and $M\approx2.61$; $P(1,2)<0$ and $P(2,2)<0$,
  \item $P(N-1,N)=0$ for $N\approx1.97$ and $N\approx20.56$; $P(1,2)<0$ and $P(15,16)<0$.
\end{itemize}
This shows that the desired inequality is true on the boundary and therefore everywhere.

\section{Case: $1\leq N\leq 2$ and $(N+1)/2\leq M\leq 2$}
In the next 2 sections we will have to deal with cases for which an equilateral triangle ($N=M=1$) is
one of the possible triangles. We take $\eps=1$, since only the eigenfunction of the equilateral 
triangle can give the constant $9$ in the Theorem \ref{main}. We also need for all the other coefficients
to vanish near the equilateral triangle. Let us take $\gamma=0$ and $\alpha=\beta$. 
We just have to choose the common value for $\alpha$ and $\beta$. Due to the nature of the already very complicated 
calculations we cannot afford to pick a very complicated coefficient, thus we take $\alpha=(N+M-2)/2$. 
This choice has one additional advantage.
In this case we are working with the eigenfunctions of the following triangles: One 
equilateral and two right triangles with shortest side $(0,0)-(1,0)$. Hence we have a symmetry
about $a=1/2$, or in terms of $M$ and $N$, about $M=N$. Therefore it is natural to introduce the rotated  
coordinates $U=(M+N)/2-1$ and $V=(N-M)/2$. Note also that $\alpha=\beta=U$.
This also moves the equilateral triangle to the origin.

After applying these transformations the inequality (\ref{horrible}) becomes
\begin{small}
\begin{gather}
  \begin{split}
    0\geq&P(U,V)=U^2\Bigl( 3293385188722144 - 451048860827136{\sqrt{3}} - 
     952832984463360\pi  
\\&\;- 874782993324240{\sqrt{3}}\pi  + 
     463182700780800{\pi }^2 - 4817666363010084U 
\\&\;+ 
     916192998555120{\sqrt{3}}U + 710514087323040{\sqrt{3}}\pi U - 
     330844786272000{\pi }^2U 
\\&\;-1072431834636645U^2 + 
     346350633108480{\sqrt{3}}\pi U^2 - 33084478627200{\pi }^2U^2
     \Bigr)  
\\+& 9V^2\Bigl( 44112638169600{\pi }^2 + 
     355514206276944{\sqrt{3}}U + 105870331607040\pi U 
\\&\;+ 
     177818984344461U^2 + 36504201333600{\sqrt{3}}\pi U^2 + 
     25732372265600{\pi }^2U^2 \Bigr)
  \end{split}
\end{gather}
\end{small}
This is  a polynomial of degree $4$ in $U$ and of degree $2$ in $V$. Hence we expect to be able to 
solve $\partial_V P(U,V)=0$. (In fact $\partial_V P$ is equal to $V$ times irreducible 
quadratic polynomial in $U$.) Therefore we have exactly one solution  
$V=0$, or $N=M$. But this is a boundary of the region, so we only have to check 
the boundary values. 

This time the boundary conditions are:  $M=N$, $M=(N+1)/2$ and $N=2$. After
changing variables to $U$ and $V$ these  become $V=0$, $U=3V$ and $U+V=1$, respectively.
Each time we get a polynomial of degree $4$. Thus we proceed as in the previous section. 
\begin{itemize}
  \item $P(U,0)=0$ for $U=0$ (double root), $U\approx5.65$ and $U\approx-0.24$; $P(0,0)=0$ and $P(1,0)<0$,
  \item $P(3V,V)=0$ for $V=0$ (double root), $V\approx0.55$ and $V\approx-0.04$; $P(0,0)=0$ and $P(3/4,1/4)<0$,
  \item $P(1-V,V)=0$ for $V\approx-0.52$ and $V\approx0.29$ ($2$ complex roots); $P(1,0)<0$ and $P(3/4,1/4)<0$.
\end{itemize}
Hence the inequality is true on the boundary and so also inside of the region.

\section{Case: $1\leq N\leq 2$ and $1\leq M\leq (N+1)/2$}
Here we take $\eps=1$, $\beta=0$ and $\alpha=\gamma=(N+M-2)/\sqrt2$.
Even though the symmetry described in the previous section does not exist here,
we will still use the same rotated coordinates $U=(M+N)/2-1$ and $V=(N-M)/2$.
This time the inequality (\ref{horrible}) becomes:
\begin{small}
\begin{gather}
  \begin{split}
0\geq&P(U,V)=32133332U^2\Bigl( -1898955433 - 549628092{\sqrt{6}} + 
     103783680{\sqrt{2}}\pi 
\\&\;\;\;- 22702680{\sqrt{3}}\pi  + 
     345945600{\pi }^2 - 1063944882U + 222614730{\sqrt{6}}U 
\\&\;\;\;+ 
     259459200{\sqrt{2}}\pi U- 136216080{\sqrt{3}}\pi U + 
     115315200{\pi }^2U - 531972441U^2 
\\&\;\;\;+ 
     113513400{\sqrt{3}}\pi U^2+ 172972800{\pi }^2U^2 \Bigr)  
\\&- 
  64266664\Bigl( 824442138{\sqrt{6}} - 155675520{\sqrt{2}}\pi  - 
     2201993543U 
\\&\;\;\;+ 158918760{\sqrt{3}}\pi U + 403603200{\pi }^2U
     \Bigr) \left( U + U^2 \right) V 
\\&+ 
  3759599844\Bigl( 1478400{\pi }^2 + 10405746{\sqrt{6}}U + 
     5765760{\sqrt{2}}\pi U - 4546773U^2 
\\&\;\;\;+ 
     4074840{\sqrt{3}}\pi U^2 + 3449600{\pi }^2U^2 \Bigr) V^2
  \end{split}
\end{gather}
\end{small}

Note that this is still a polynomial of degree $2$ in $V$ and  therefore we proceed as in the
previous section. Unfortunately this time the only solution of $\partial_V P=0$ is $V$ equal
to a rational function of $U$ with an irreducible denominator of degree $2$. 
Hence, by plugging this into $\partial_U P=0$ we get a rational equation with squared 
irreducible polynomial of degree $2$ in the denominator. Thus, this equation is 
equivalent to the numerator being $0$. Fortunately the numerator is a solvable polynomial
of degree $7$ with $4$ imaginary roots and $3$ real roots ($0$, $\approx-0.18$ and $\approx1.8$).

Here  we have the following bounds $U=3V$ (equivalent to $M=(N+1)/2$),
$U+V=1$ ($N=2$), and $U=V$ ($M=1$). So this is a triangle with vertices
$(0,0)$, $(3/4,1/4)$ and $(1/2,1/2)$. Hence neither critical point is inside of this region.
Therefore we have to check the boundary values and like before this means we have to find
the roots of certain polynomials of degree $4$ and values at the endpoints.
\begin{itemize}
  \item $P(V,V)=0$ for $V=0$ (double root), $V\approx-0.27$ and $V\approx0.64$; 
    $P(0,0)=0$ and $P(1/2,1/2)<0$,
  \item $P(3V,V)=0$ for $V=0$ (double root), $V\approx-0.06$ and $V\approx0.51$;
    $P(0,0)=0$ and $P(3/4,1/4)<0$,
  \item $P(U,1-U)=0$ for $U\approx0.48$ and $U\approx0.79$ ($2$ complex roots); 
    $P(1/2,1/2)<0$ and $P(3/4,1/4)<0$.
\end{itemize}
Thus inequality is true. 

\section{Case: $M\geq 15$}
For this case  we will use a method different than all other cases. Since we are dealing 
with the triangles for which two sides are long and almost equal, we will estimate the  
eigenvalue by the eigenvalue of a circular sector contained in the triangle $T$.

Let us denote the angle between the sides of length $N$ and $M$ by $\gamma$. First we 
take the isosceles triangle with angle $\gamma$ between the sides of length $M$. We can certainly
put this triangle inside the triangle $T$. Since the shortest side of this isosceles triangle has length no larger 
than $1$, the altitude $h$ satisfies $h\geq\sqrt{M^2-1/4}$.

Let us denote a circular sector with angle $\alpha$ and radius $r$ by $S(\alpha,r)$. 
It is known (see \cite{PS}), that the first eigenvalue of the sector $S(\alpha,r)$ is 
$j^2_{\pi/\alpha}r^{-2}$, where $j_{\nu}$ is the first zero of the Bessel function
$J_{\nu}(x)$ of order $\nu$. 

It is clear that we can put a sector $S(\gamma,h)$ inside the triangle $T$. Hence, by domain
monotonicity we have
\begin{gather}
  \lambda_T\leq \lambda_{S(\gamma,h)}={j^2_{\pi/\gamma}\over h^2}.
\end{gather}
We need to prove that
\begin{gather}
  {j^2_{\pi/\gamma}\over h^2}\leq {\pi^2 L_T^2\over 9A_T^2}.
\end{gather}

We have $L_T=1+M+N\geq 2N$ and $A_T= \sin(\gamma) NM/2\leq \gamma NM/2$.
Therefore, it is enough to prove that
\begin{gather}
  {9j^2_{\pi/\gamma}(\gamma NM/2)^2\over (M^2-1/4)^2 (2N)^2}\leq 1,
\end{gather}
or that
\begin{gather}\label{est}
  {9j^2_{\pi/\gamma}\gamma^2 M^2\over 16\pi^2(M^2-1/4)^2}\leq 1.
\end{gather}

To find the bound for $j_{\nu}$ we will use the estimate obtained in \cite{QW}
\begin{gather}
  j_{\nu}\leq \nu-{a_1\over \sqrt[3]{2}}\nu^{1/3}+{3a_1^2\sqrt[3]{2}\over 20}\nu^{-1/3},
\end{gather}
where $a_1\approx-2.338$ is the first negative zero of the Airy function.
Hence we have
\begin{gather}
  {j_{\nu}\over \nu}\leq 1+2\nu^{-2/3}+2\nu^{-4/3}.
\end{gather}

Therefore
\begin{gather}
  {9j^2_{\pi/\gamma}\gamma^2 M^2\over 16\pi^2(M^2-1/4)^2}\leq
  \left(1+2\left(\frac\gamma\pi\right)^{2/3}+2\left(\frac\gamma\pi\right)^{4/3}\right)^2
  {9M^2\over 16(M^2-1/4)}.
\end{gather}

This last expression is increasing with $\gamma$ as can be easily verified by differentiating. 
Given $M$, the angle $\gamma$ is maximized
for the isosceles triangle, hence $\gamma\leq 2\sin^{-1}(1/2M)$.
In order to arrive at (\ref{est}), it is enough to show 
\begin{gather}
  \left(1+2\left(\frac{2\sin^{-1}(1/2M)}\pi\right)^{2/3}+2\left(\frac{2\sin^{-1}(1/2M)}
  \pi\right)^{4/3}\right)^2\!\!
  {9M^2\over 16(M^2-1/4)}\leq \! 1.
\end{gather}
 It is easy to check that the function on the left side is decreasing with $M$, 
 and that for $M=15$ inequality is true. Hence this is true for any triangle with 
 $M\geq15$.

 Note also that if $M\To\infty$, then the whole expression tends to $9/16$. This shows that the constant $16$
 in the lower bound in Theorem \ref{main} is optimal.
 
\section{Script in Mathematica}
Here we give the script written in Mathematica to handle all the cumbersome calculations
included in Sections 2 to 5. It is important to note that all the calculations are done
symbolically. Only the exact values of the roots of all the polynomials are at  end
converted to numerical form.
\begin{small}
\begin{verbatim}
(* Section 2 *)
(* isosceles triangle with vertices (0,0), (1,0) and (Sqrt[3],0) *)
g[x_,y_]=Sin [Sqrt[3]\[Pi] y]Sin[\[Pi] x/3] +  \
        Sin[\[Pi] y/Sqrt[3]]Sin[5\[Pi] x/3] +  \
        Sin[2\[Pi] y/Sqrt[3]]Sin[4\[Pi] x/3];
(* other right triangles *)
g2[x_,y_]=g[1-x,y];          
g3[x_,y_]=g[Sqrt[3]y,Sqrt[3]x];
(* test functions obtained from right triangles *)
\[CurlyPhi]1=g[x-(a y /b),Sqrt[3]y/b];
\[CurlyPhi]2=g2[x-((a-1) y /b),Sqrt[3]y/b];
\[CurlyPhi]3=g3[x-(a y /b),y/(Sqrt[3]b)];
(* equilateral triangle after linear transformation *)
\[Phi]:=Sin[2\[Pi]y/b]-Sin[2\[Pi](x+(1-a)y/b)]+Sin[2\[Pi](x-a y/b)]; 
(* final test function *)
\[Psi]=\[Alpha] \[CurlyPhi]1 + \[Beta] \[CurlyPhi]2 +  \
        \[Gamma] \[CurlyPhi]3 + \[Epsilon] \[Phi];
grad=Simplify[Integrate[D[\[Psi],x]^2+D[\[Psi],y]^2,{y,0,b},  \
        {x,a y/b, (a-1) y/b+1}]];
int=Simplify[ Integrate[\[Psi]^2,{y,0,b},{x,a y/b , (a-1)y/b +1}]];
(* we have to prove that this is <= 0 *)
in=9b^2grad-4\[Pi]^2(1+Sqrt[a^2+b^2]+Sqrt[(a-1)^2+b^2])^2int;
(* change from (a, b) to (M, N) and cancel b *)
in2=Simplify[in/b /. b^2 -> M^2 - a^2 /. a -> (M^2 - N^2 + 1)/2,  \
        (N > 0) && (M > 0)];
(* inequality (2.9) *)
Simplify[308788467187200in/b]

(* Section 3 *)
W=in2/. \[Epsilon] -> 0 /. \[Gamma] -> -1/6 /. \[Beta] -> 0 /.  \
        \[Alpha] -> 1;
(* Inequality (3.1) *)
Apart[1383782400W]
(* Critical point *)
Reduce[(D[W, M] == 0) && (D[W, N] == 0), {M, N}] // N
(* Boundary : roots and endpoints *)
Reduce[W == 0 /. N -> 2] // N
Reduce[W == 0 /. M -> N - 1] // N
Reduce[W == 0 /. M -> N] // N
Reduce[W == 0 /. M -> 15] // N
W /. M -> {1, 2} /. N -> 2 // N
W /. M -> 15 /. N -> {15, 16} // N

(* Section 4 *)
W=in2/. \[Epsilon] -> 1 /. \[Gamma] -> 0 /. \[Beta] -> \[Alpha] /.\
        \[Alpha] -> (N + M - 2)/2;
pol = W/. M -> U - V /. N -> U + V /. U -> U + 1;
(* inequality (4.1) *)
Apart[22056319084800pol, V]
(* Critical point *) 
Reduce[D[pol, V] == 0, V] // N
(* Boundary : roots and endpoints *)
Reduce[pol == 0 /. V -> 0] // N
Reduce[pol == 0 /. U -> 1 - V] // N
Reduce[pol == 0 /. U -> 3V] // N
pol /. V -> 0 /. U -> {0, 1} // N
pol /. V -> 1/4 /. U -> 3/4 // N

(* Section 5 *)
W=in2/. \[Epsilon] -> 1 /. \[Beta] -> 0 /. \[Gamma] -> \[Alpha] /.\
        \[Alpha] -> (M + N - 2)/Sqrt[2];
pol = W /. M -> U - V /. N -> U + V /. U -> U + 1;
(* inequality (5.1) *)
Apart[9609600pol, V]
(* Critical points *)
Vs = Solve[D[pol, V] == 0, V];
Reduce[D[pol, V] == 0, V, Reals] 
(* denominator with complx roots only*)
Reduce[Denominator[Together[D[pol, U] /. Vs]] == 0] // N 
(* polynomial of degree 7 in U *)
Reduce[Numerator[Together[
      D[pol, U] /. Vs]] == 0] // N
(* Boundary : roots and endpoints *)
Reduce[pol == 0 /. U -> 3V] // N
Reduce[pol == 0 /. U -> V] // N
Reduce[pol == 0 /. V -> 1 - U] // N
pol /. U -> 1 - V /. V -> {1/4, 1/2} // N
pol /. U -> 0 /. V -> 0 // N
\end{verbatim}
\end{small}

\section*{Acknowledgements}
The author wants to thank his thesis advisor, Professor Rodrigo Ba\~nuelos, for
the support and quidance on this paper, which is a part of author's Ph. D. thesis. 

\bibliographystyle{amsplain}

\begin{thebibliography}{10}
  \bibitem{AF} P. Antunes and P. Freitas; \textit{New bounds for the principal Dirichlet
    eigenvalue of planar regions}, preprint.
  \bibitem{F} P. Freitas; \textit{Upper and lower bounds for the first Dirichlet 
    eigenvalue of a triangle}, Proc. Amer. Math. Soc., posted on January 6, 2006, 
    PII S 0002-9939(06)08339-0(to appear in print).
  \bibitem{M} E. Makai; \textit{On the principal frequency of a membrane and the
    torsion rigidity of a beam}, pp. 227-231 in Studies in mathematical analysis
    and related topics, Essays in honor of George P\'olya, Stanford Univ. Press, Stanford 1962.
  \bibitem{Mc} B. McCartin; \textit{Eigenstructure of the equilateral triangle. I. The Dirichlet problem.},
    SIAM Rev. {\bf45} (2003), {\it no. 2}, 267--287.
  \bibitem{O} R. Osserman; \textit{A note on Hayman's theorem on the bass note of a drum},
    Comment. Math. Helvetici {\bf52} (1977), 545--555.
  \bibitem{Po} G. P\'olya; \textit{Two more inequalities between physical and geometrical
    quantities}, J. Indian Math. Soc. (N.S.){\bf24} (1960), 413--419.
  \bibitem{PS} G. P\'olya and G. Szeg\"o; \textit{Isoperimetric inequalities in mathematical
    physics}, Annals of Mathematical Studies {\bf27}, Princeton University Press, 1951.
  \bibitem{QW} C. K. Qu and R. Wong; \textit{``Best possible'' upper and lower bounds for the zeros of the Bessel function 
    $J_v(x)$}, Trans. Amer. Math. Soc. {\bf351} (1999), 2833--2859.
\end{thebibliography}

\end{document}